\documentclass[10.9pt,twoside]{amsart}
\usepackage{amsmath, amsthm, amscd, amsfonts, amssymb, graphicx, color}
\usepackage[bookmarksnumbered, colorlinks, plainpages]{hyperref}

\theoremstyle{definition}

\theoremstyle{remark}

\numberwithin{equation}{section}

\begin{document}
\setcounter{page}{1}
\begin{center}
{\bf  THE ARENS REGULARITY AND WEAK \\TOPOLOGICAL CENTER OF MODULE ACTIONS}
\end{center}

\title[]{}
\author[]{KAZEM HAGHNEJAD AZAR}

\address{}

\dedicatory{}



\begin{abstract} Let $A$ be a Banach algebra and $A^{**}$ be the second dual of it. We define $\tilde{Z}_1(A^{**})$ as a weak topological center of $A^{**}$ with respect to first Arens product and find some relations between it and the topological center of $A^{**}$. We also extend this new definition  into the module actions and find relationship between weak topological center of module actions and reflexivity or Arens regularity of some Banach algebras, and we investigate some applications of this new definition in the weak amenability of some Banach algebras.

\end{abstract} \maketitle

\section{\bf  preliminaries and
Introduction }

\noindent  As is well-known [1], the second dual $A^{**}$ of $A$ endowed with the either Arens multiplications is a Banach algebra. The constructions of the two Arens multiplications in $A^{**}$ lead us to definition of topological centers for $A^{**}$ with respect both Arens multiplications. The topological centers of Banach algebras, module actions and applications of them  were introduced and discussed in [6, 8, 13, 14, 15, 16, 17, 21, 23, 24], and they have attracted by some attentions. In this paper we will introduce a new concept as a weak topological center of Banach algebras, module actions and  investigate its relationship  with topological centers of Banach algebras, module actions and we will obtain some useful conclusions in Banach algebras with some example in the dual groups. In final section, we give some applications of weak topological center of Banach algebras in  weak amenability of Banach algebras.\\
Now we introduce some notations and definitions that we used
throughout  this paper. Let $A$ be  a Banach algebra and $A^*$,
$A^{**}$, respectively, are the first and second dual of $A$.  For $a\in A$
 and $a^\prime\in A^*$, we denote by $a^\prime a$
 and $a a^\prime$ respectively, the functionals on $A^*$ defined by $<a^\prime a,b>=<a^\prime,ab>=a^\prime(ab)$ and $<a a^\prime,b>=<a^\prime,ba>=a^\prime(ba)$ for all $b\in A$.
   The Banach algebra $A$ is embedded in its second dual via the identification
 $<a,a^\prime>$ - $<a^\prime,a>$ for every $a\in
A$ and $a^\prime\in
A^*$.
 \noindent Let $X,Y,Z$ be normed spaces and $m:X\times Y\rightarrow Z$ be a bounded bilinear mapping. Arens in [1] offers two natural extensions $m^{***}$ and $m^{t***t}$ of $m$ from $X^{**}\times Y^{**}$ into $Z^{**}$ as following\\
1. $m^*:Z^*\times X\rightarrow Y^*$,~~~~~given by~~~$<m^*(z^\prime,x),y>=<z^\prime, m(x,y)>$ ~where $x\in X$, $y\in Y$, $z^\prime\in Z^*$,\\
 2. $m^{**}:Y^{**}\times Z^{*}\rightarrow X^*$,~~given by $<m^{**}(y^{\prime\prime},z^\prime),x>=<y^{\prime\prime},m^*(z^\prime,x)>$ ~where $x\in X$, $y^{\prime\prime}\in Y^{**}$, $z^\prime\in Z^*$,\\
3. $m^{***}:X^{**}\times Y^{**}\rightarrow Z^{**}$,~ given by~ ~ ~$<m^{***}(x^{\prime\prime},y^{\prime\prime}),z^\prime>$  $=<x^{\prime\prime},m^{**}(y^{\prime\prime},z^\prime)>$\\ ~where ~$x^{\prime\prime}\in X^{**}$, $y^{\prime\prime}\in Y^{**}$, $z^\prime\in Z^*$.\\
The mapping $m^{***}$ is the unique extension of $m$ such that $x^{\prime\prime}\rightarrow m^{***}(x^{\prime\prime},y^{\prime\prime})$ from $X^{**}$ into $Z^{**}$ is $weak^*-weak^*$ continuous for every $y^{\prime\prime}\in Y^{**}$, but the mapping $y^{\prime\prime}\rightarrow m^{***}(x^{\prime\prime},y^{\prime\prime})$ is not in general $weak^*-weak^*$ continuous from $Y^{**}$ into $Z^{**}$ unless $x^{\prime\prime}\in X$. Hence the first topological center of $m$ may  be defined as following
$$Z_1(m)=\{x^{\prime\prime}\in X^{**}:~~y^{\prime\prime}\rightarrow m^{***}(x^{\prime\prime},y^{\prime\prime})~~is~~weak^*-weak^*-continuous\}.$$
Let now $m^t:Y\times X\rightarrow Z$ be the transpose of $m$ defined by $m^t(y,x)=m(x,y)$ for every $x\in X$ and $y\in Y$. Then $m^t$ is a continuous bilinear map from $Y\times X$ to $Z$, and so it may be extended as above to $m^{t***}:Y^{**}\times X^{**}\rightarrow Z^{**}$.
 The mapping $m^{t***t}:X^{**}\times Y^{**}\rightarrow Z^{**}$ in general is not equal to $m^{***}$, see [1], if $m^{***}=m^{t***t}$, then $m$ is called Arens regular. The mapping $y^{\prime\prime}\rightarrow m^{t***t}(x^{\prime\prime},y^{\prime\prime})$ is $weak^*-weak^*$ continuous for every $y^{\prime\prime}\in Y^{**}$, but the mapping $x^{\prime\prime}\rightarrow m^{t***t}(x^{\prime\prime},y^{\prime\prime})$ from $X^{**}$ into $Z^{**}$ is not in general  $weak^*-weak^*$ continuous for every $y^{\prime\prime}\in Y^{**}$. So we define the second topological center of $m$ as
$$Z_2(m)=\{y^{\prime\prime}\in Y^{**}:~~x^{\prime\prime}\rightarrow m^{t***t}(x^{\prime\prime},y^{\prime\prime})~~is~~weak^*-weak^*-continuous\}.$$
It is clear that $m$ is Arens regular if and only if $Z_1(m)=X^{**}$ or $Z_2(m)=Y^{**}$. Arens regularity of $m$ is equivalent to the following
$$\lim_i\lim_j<z^\prime,m(x_i,y_j)>=\lim_j\lim_i<z^\prime,m(x_i,y_j)>,$$
whenever both limits exist for all bounded sequences $(x_i)_i\subseteq X$ , $(y_i)_i\subseteq Y$ and $z^\prime\in Z^*$, see [6, 14, 18].\\
 The regularity of a normed algebra $A$ is defined to be the regularity of its algebra multiplication when considered as a bilinear mapping. Let $a^{\prime\prime}$ and $b^{\prime\prime}$ be elements of $A^{**}$, the second dual of $A$. By $Goldstin^,s$ Theorem [4, P.424-425], there are nets $(a_{\alpha})_{\alpha}$ and $(b_{\beta})_{\beta}$ in $A$ such that $a^{\prime\prime}=weak^*-\lim_{\alpha}a_{\alpha}$ ~and~  $b^{\prime\prime}=weak^*-\lim_{\beta}b_{\beta}$. So it is easy to see that for all $a^\prime\in A^*$,
$$\lim_{\alpha}\lim_{\beta}<a^\prime,m(a_{\alpha},b_{\beta})>=<a^{\prime\prime}b^{\prime\prime},a^\prime>$$ and
$$\lim_{\beta}\lim_{\alpha}<a^\prime,m(a_{\alpha},b_{\beta})>=<a^{\prime\prime}ob^{\prime\prime},a^\prime>,$$
where $a^{\prime\prime}.b^{\prime\prime}$ and $a^{\prime\prime}ob^{\prime\prime}$ are the first and second Arens products of $A^{**}$, respectively, see [6, 14, 18].\\
The mapping $m$ is left strongly Arens irregular if $Z_1(m)=X$ and $m$ is right strongly Arens irregular if $Z_2(m)=Y$.\\
Regarding $A$ as a Banach $A-bimodule$, the operation $\pi:A\times A\rightarrow A$ extends to $\pi^{***}$ and $\pi^{t***t}$ defined on $A^{**}\times A^{**}$. These extensions are known, respectively, as the first (left) and the second (right) Arens products, and with each of them, the second dual space $A^{**}$ becomes a Banach algebra. In this situation, we shall also simplify our notations. So the first (left) Arens product of $a^{\prime\prime},b^{\prime\prime}\in A^{**}$ shall be simply indicated by $a^{\prime\prime}b^{\prime\prime}$ and defined by the three steps:
 $$<a^\prime a,b>=<a^\prime ,ab>,$$
  $$<a^{\prime\prime} a^\prime,a>=<a^{\prime\prime}, a^\prime a>,$$
  $$<a^{\prime\prime}b^{\prime\prime},a^\prime>=<a^{\prime\prime},b^{\prime\prime}a^\prime>.$$
 for every $a,b\in A$ and $a^\prime\in A^*$. Similarly, the second (right) Arens product of $a^{\prime\prime},b^{\prime\prime}\in A^{**}$ shall be  indicated by $a^{\prime\prime}ob^{\prime\prime}$ and defined by :
 $$<a oa^\prime ,b>=<a^\prime ,ba>,$$
  $$<a^\prime oa^{\prime\prime} ,a>=<a^{\prime\prime},a oa^\prime >,$$
  $$<a^{\prime\prime}ob^{\prime\prime},a^\prime>=<b^{\prime\prime},a^\prime ob^{\prime\prime}>.$$
  for all $a,b\in A$ and $a^\prime\in A^*$.\\
  The regularity of a normed algebra $A$ is defined to be the regularity of its algebra multiplication when considered as a bilinear mapping. Let $a^{\prime\prime}$ and $b^{\prime\prime}$ be elements of $A^{**}$, the second dual of $A$. By $Goldstine^,s$ Theorem [4, P.424-425], there are nets $(a_{\alpha})_{\alpha}$ and $(b_{\beta})_{\beta}$ in $A$ such that $a^{\prime\prime}=weak^*-\lim_{\alpha}a_{\alpha}$ ~and~  $b^{\prime\prime}=weak^*-\lim_{\beta}b_{\beta}$. So it is easy to see that for all $a^\prime\in A^*$,
$$\lim_{\alpha}\lim_{\beta}<a^\prime,\pi(a_{\alpha},b_{\beta})>=<a^{\prime\prime}b^{\prime\prime},a^\prime>$$ and
$$\lim_{\beta}\lim_{\alpha}<a^\prime,\pi(a_{\alpha},b_{\beta})>=<a^{\prime\prime}ob^{\prime\prime},a^\prime>,$$
where $a^{\prime\prime}b^{\prime\prime}$ and $a^{\prime\prime}ob^{\prime\prime}$ are the first and second Arens products of $A^{**}$, respectively, see [6, 14, 18].\\
We find the usual first and second topological center of $A^{**}$, which are
  $$Z_1(A^{**})=Z^\ell_1(A^{**})=\{a^{\prime\prime}\in A^{**}: b^{\prime\prime}\rightarrow a^{\prime\prime}b^{\prime\prime}~ is~weak^*-weak^*~continuous\},$$
   $$Z_2(A^{**})=Z_2^r(A^{**})=\{a^{\prime\prime}\in A^{**}: a^{\prime\prime}\rightarrow a^{\prime\prime}ob^{\prime\prime}~ is~weak^*-weak^*~continuous\}.$$\\
 An element $e^{\prime\prime}$ of $A^{**}$ is said to be a mixed unit if $e^{\prime\prime}$ is a
right unit for the first Arens multiplication and a left unit for
the second Arens multiplication. That is, $e^{\prime\prime}$ is a mixed unit if
and only if,
for each $a^{\prime\prime}\in A^{**}$, $a^{\prime\prime}e^{\prime\prime}=e^{\prime\prime}o a^{\prime\prime}=a^{\prime\prime}$. By
[4, p.146], an element $e^{\prime\prime}$ of $A^{**}$  is  mixed
      unit if and only if it is a $weak^*$ cluster point of some BAI $(e_\alpha)_{\alpha \in I}$  in
      $A$.\\
Let now $B$ be a Banach $A-bimodule$, and let
$$\pi_\ell:~A\times B\rightarrow B~~~and~~~\pi_r:~B\times A\rightarrow B.$$
be the right and left module actions of $A$ on $B$. Then $B^{**}$ is a Banach $A^{**}-bimodule$ with module actions
$$\pi_\ell^{***}:~A^{**}\times B^{**}\rightarrow B^{**}~~~and~~~\pi_r^{***}:~B^{**}\times A^{**}\rightarrow B^{**}.$$
Similarly, $B^{**}$ is a Banach $A^{**}-bimodule$ with module actions\\
$$\pi_\ell^{t***t}:~A^{**}\times B^{**}\rightarrow B^{**}~~~and~~~\pi_r^{t***t}:~B^{**}\times A^{**}\rightarrow B^{**}.$$
Assume that  $B$ is a Banach $A-bimodule$. We say that  $B$ factors on the left (right) with respect to $A$ if $B=BA~(B=AB)$. We say that $B$ factors on both sides with respect to $A$, if $B=BA=AB$.\\
Let $B$ be a left Banach $A-module$ and  $e$ be a left  unit element of $A$. Then we say that $e$ is a left unit (resp. weakly left unit)  $A-module$ for $B$, if $\pi_\ell(e,b)=b$ (resp. $<b^\prime , \pi_\ell(e,b)>=<b^\prime , b>$ for all $b^\prime\in B^*$) where $b\in B$. The definition of right unit (resp. weakly right unit) $A-module$ is similar.\\
We say that a Banach $A-bimodule$ $B$ is a unital $A-module$, if $B$ has left and right unit $A-module$ that are equal then we say that $B$ is unital $A-module$.\\
Let $B$ be a left Banach $A-module$ and $(e_{\alpha})_{\alpha}\subseteq A$ be a LAI [resp. weakly left approximate identity(=WLAI)] for $A$. We say that $(e_{\alpha})_{\alpha}$ is left approximate identity $A-module$ ($=LAI~A-module$)[ resp. weakly left approximate identity $A-module$ (=$WLAI~A-module$)] for $B$, if for all $b\in B$, we have $\pi_\ell (e_{\alpha},b) \stackrel{} {\rightarrow}
b$ ( resp. $\pi_\ell (e_{\alpha},b) \stackrel{w} {\rightarrow}
b$). The definition of the right approximate identity ($=RAI~A-module$)[ resp. weakly right approximate identity ($=WRAI~A-module$)] is similar.\\
We say that $(e_{\alpha})_{\alpha}$ is a approximate identity $A-module$ ($=AI~A-module$)[ resp. weakly approximate identity $A-module$ ($WAI~A-module$)] for $B$, if $B$ has left and right approximate identity $A-module$ [ resp. weakly left and right approximate identity $A-module$] that are equal.\\
Let $(e_{\alpha})_{\alpha}\subseteq A$ be $weak^*$ left approximate identity for $A^{**}$. We say that $(e_{\alpha})_{\alpha}$ is $weak^*$ left approximate identity $A^{**}-module$ ($=W^*LAI~A^{**}-module$) for $B^{**}$, if for all $b^{\prime\prime}\in B^{**}$, we have $\pi_\ell^{***} (e_{\alpha},b^{\prime\prime}) \stackrel{w^*} {\rightarrow}
b^{\prime\prime}$. The definition of the $weak^*$ right approximate identity $A^{**}-module$($=W^*RAI~A^{**}-module$) is similar.\\
 We say that $(e_{\alpha})_{\alpha}$ is a $weak^*$ approximate identity $A^{**}-module$ ($=W^*AI~A^{**}-module$) for $B^{**}$, if $B^{**}$ has $weak^*$ left and right approximate identity $A^{**}-module$  that are equal.\\
A functional $a^\prime$ in $A^*$ is said to be $wap$ (weakly almost
 periodic) on $A$ if the mapping $a\rightarrow a^\prime a$ from $A$ into
 $A^{*}$ is weakly compact. Pym in [18] showed that  this definition to the equivalent following condition\\
 For any two net $(a_{\alpha})_{\alpha}$ and $(b_{\beta})_{\beta}$
 in $\{a\in A:~\parallel a\parallel\leq 1\}$, we have\\
$$\\lim_{\alpha}\\lim_{\beta}<a^\prime,a_{\alpha}b_{\beta}>=\\lim_{\beta}\\lim_{\alpha}<a^\prime,a_{\alpha}b_{\beta}>,$$
whenever both iterated limits exist. The collection of all $wap$
functionals on $A$ is denoted by $wap(A)$. Also we have
$a^{\prime}\in wap(A)$ if and only if $<a^{\prime\prime}b^{\prime\prime},a^\prime>=<a^{\prime\prime}ob^{\prime\prime},a^\prime>$ for every $a^{\prime\prime},~b^{\prime\prime} \in
A^{**}$. \\
\noindent Let $B$ be a left Banach $A-module$. Then, $b^\prime\in B^*$ is said to be left weakly almost periodic functional if the set $\{\pi_\ell(b^\prime,a):~a\in A,~\parallel a\parallel\leq 1\}$ is relatively weakly compact. We denote by $wap_\ell(B)$ the closed subspace of $B^*$ consisting of all the left weakly almost periodic functionals in $B^*$.\\
The definition of the right weakly almost periodic functional ($=wap_r(B)$) is the same.\\
By [6, 14, 18], the definition of $wap_\ell(B)$ is equivalent to the following $$<\pi_\ell^{***}(a^{\prime\prime},b^{\prime\prime}),b^\prime>=
<\pi_\ell^{t***t}(a^{\prime\prime},b^{\prime\prime}),b^\prime>$$
for all $a^{\prime\prime}\in A^{**}$ and $b^{\prime\prime}\in B^{**}$.
Thus, we can write \\
$$wap_\ell(B)=\{ b^\prime\in B^*:~<\pi_\ell^{***}(a^{\prime\prime},b^{\prime\prime}),b^\prime>=
<\pi_\ell^{t***t}(a^{\prime\prime},b^{\prime\prime}),b^\prime>~~$$$$for~~all~~a^{\prime\prime}\in A^{**},~b^{\prime\prime}\in B^{**}\}.$$\\
In many parts of this paper, for left or right Banach $A-module$, we take $\pi_\ell(a,b)=ab$ and  $\pi_r(b,a)=ba$, for every $a\in A$ and $b\in B$.\\ This paper is organized as follows:\\
{\bf a}) Suppose that $A$ is a Banach algebra and $A^{***}A^{**}\subseteq A^*$. Then the topological centers of $A^{**}$ (=${{Z}_1}(A^{**})$)   and weak topological centers of $A^{**}$ (=$\tilde{{Z}_1^\ell}(A^{**})$) are $A^{**}$.\\
{\bf b}) Let $A$ be a Banach algebra and let the left weak topological center of $A^{**}$ be $A$. Then $A$ is a right ideal in $A^{**}$.\\
{\bf c}) Let $B$ a Banach $A-bimodule$. Then for $n\geq 2$, we have the following assertions.\\
i) ~If $B^{(n)}B^{(n-1)}\subseteq A^{(n-2)}$, then ${Z}^\ell_{A^{(n-1)}}(B^{(n-1)})={\tilde{{Z}}^\ell}_{A^{(n-1)}}(B^{(n-1)})$.\\
ii) If $B^{(n)}A^{(n-1)}\subseteq B^{(n-2)}$ and $B^{(n-1)}$ has a left unit  $A^{(n-1)}-module$, then $B$ is reflexive.\\
{\bf d}) Let $A$ be a Banach algebra. Then\\
i) If $A^{***}A^{**}\subseteq \widetilde{wap}(A)$, then $A^{**}A^{**}\subseteq {{Z}_1}(A^{**})$ where
$$\widetilde{wap}(A)=\{a^{\prime\prime\prime}\in A^{***}:~a^{\prime\prime\prime}a^{\prime\prime}:A^{**}\rightarrow
\mathbb{C}~is ~weak^*~continuous~for~all~a^{\prime\prime}\in A^{**}\}.$$
ii) If $A^{**}=A^{**}\tilde{{Z}_1}^\ell(A^{**})$, then $A^{***}A^{**}\subseteq \widetilde{wap}(A)$.\\
{\bf e}) Let $B$ a right  Banach $A-module$ and $e^{(n)}\in A^{(n)}$ be a right unit for $B^{(n)}$  where $n\geq 2$. Then we have the following assertions.\\
i)~~~ If ${ {Z}^\ell}_{e^{(n)}}(B^{(n)})=B^{(n)}$, then $A^{(n-2)}B^{(n-1)}= B^{(n-1)}$ where
$$\tilde{ {Z}^\ell}_{e^{n}}(B^{n})=\{ b^{n}\in B^{n}:~e^{\prime\prime}\rightarrow b^{n}e^{n}~is~weak^*-to-weak~continuous\},$$
ii) If $\tilde{ {Z}^\ell}_{e^{(n)}}(B^{(n)})=B^{(n)}$, then $A^{(n-2)}B^{(n-1)}= B^{(n-1)}$ and $weak^*$ closure of \\ $A^{(n-2)}B^{(n+1)}$~ is~ $B^{(n+1)}$.\\
\noindent {\bf f})  Assume that $A$ is a Banach algebra and $\tilde{{Z}_1^\ell}(A^{(n)})=A^{(n)}$ where $n>1$. If $A^{(n)}$ is weakly amenable, then $A^{(n-2)}$ is weakly amenable.\\
 {\bf g}) Let $B$ a   Banach $A-bimodule$ and $D:A^{(n)}\rightarrow B^{(n+1)}$ be a derivation where $n\geq 0$. If $\tilde{Z}_{A^{(n)}}^\ell(B^{(n)})=B^{(n)}$, then $D^{\prime\prime}:A^{(n+2)}\rightarrow A^{(n+3)}$ is a derivation.\\
{\bf j}) Let $B$ be a   Banach $A-bimodule$ and $\tilde{{Z}^\ell}_{A^{{**}}}(B^{{**}})=B^{{**}}$. Assume that  $D:A\rightarrow B^*$ is a derivation. Then $D^{**}(A^{**})B^{**}\subseteq A^*$.\\

\begin{center}
\section{ \bf  Weak topological center of Banach algebras  }
\end{center}

\noindent In this section, we define a new concept as weak topological center of the second dual of a Banach algebra with respect to first and second Arens product and we  find its relations  with the topological centers of the second dual of a Banach algebra with respect to first and second Arens product. We establish some conditions that  when the weak topological center and topological center of the second dual of Banach algebra with respect to first Arens product coincide.\\

\noindent{\it{\bf Definition 2-1.}} Let $A$ be a Banach algebra. Then we define weak topological centers of $A^{**}$ with respect to the first and second Arens product, respectively, in the following\\
$$\tilde{ {Z}^\ell_1}(A^{**})=\{ a^{\prime\prime}\in A^{**}:~b^{\prime\prime}\rightarrow a^{\prime\prime}b^{\prime\prime}~is~weak^*-to-weak~continuous~for~all~b^{\prime\prime}\in A^{**}\},$$
$$\tilde{ {Z}^r_1}(A^{**})=\{ a^{\prime\prime}\in A^{**}:~b^{\prime\prime}\rightarrow b^{\prime\prime}a^{\prime\prime}~is~weak^*-to-weak~continuous~for~all~b^{\prime\prime}\in A^{**}\},$$
$$\tilde{ {Z}^\ell_2}(A^{**})=\{ a^{\prime\prime}\in A^{**}:~b^{\prime\prime}\rightarrow b^{\prime\prime}oa^{\prime\prime}~is~weak^*-to-weak~continuous~for~all~b^{\prime\prime}\in A^{**}\},$$
$$\tilde{ {Z}^r_2}(A^{**})=\{ a^{\prime\prime}\in A^{**}:~b^{\prime\prime}\rightarrow a^{\prime\prime}ob^{\prime\prime}~is~weak^*-to-weak~continuous~for~all~b^{\prime\prime}\in A^{**}\}.$$\\
\noindent It is clear that $\tilde{{Z}_i^\ell}(A^{**})$ and $\tilde{{Z}_i^r}(A^{**})$ are subspace of $A^{**}$ with respect to the first and second Arens products and we have also  $\tilde{{Z}_1^\ell}(A^{**})\subseteq {Z}_1(A^{**})$, ~ $\tilde{{Z}_2^r}(A^{**})\subseteq {Z}_2(A^{**})$  and if ${Z}_1(A^{**})={Z}_2(A^{**})$, it follow that $\tilde{{Z}_1^\ell}(A^{**})=\tilde{{Z}_2^\ell}(A^{**})$, $\tilde{{Z}_1^r}(A^{**})=\tilde{{Z}_2^r}(A^{**})$, and  if  ${Z}_1(A^{**})=\tilde{{Z}_1^\ell}(A^{**})=\tilde{{Z}_2^r}(A^{**})$ or $\tilde{{Z}_1^\ell}(A^{**})=\tilde{{Z}_2^r}(A^{**})={Z}_2(A^{**})$, then we conclude that ${Z}_1(A^{**})={Z}_2(A^{**})$. Obviously, if $\tilde{{Z}_1^\ell}(A^{**})= A^{**}$ or $\tilde{{Z}_1^r}(A^{**})= A^{**}$, then $A$ is Arens regular. Now let $\tilde{{Z}_1^\ell}(A^{**})={{Z}_1}(A^{**})$. Then for every subspace $B$ of $A^{**}$, we have $B\tilde{{Z}_1^\ell}(A^{**})=B{{Z}_1}(A^{**})$.\\\\
If a Banach algebra $A$ is  Arens regular or strongly Arens irregular on the left, in general, we can not conclude that $\tilde{{Z}_1^\ell}(A^{**})=A^{**}$ or  $\tilde{{Z}_1^\ell}(A^{**})=A^{}$, respectively. \\
In the following, we give  examples of some Arens regular or strongly Arens irregular Banach algebras such as $A$ that  $\tilde{{Z}_1^\ell}(A^{**})$ is equal to $A^{**}$ or no.\\
 i) Let $A$ be nonreflexive Arens regular  Banach algebra and $e^{\prime\prime}\in A^{**}$ be a left unite element of $A^{**}$. Then $\tilde{{Z}_1^\ell}(A^{**})\neq A^{**}.$\\
 ii) Suppose that  $G$ is  a locally compact group. Then $M(G)$ is left strong Arens irregular, but $\tilde{{Z}_1^\ell}{(M(G)^{**}})\neq M(G)$.\\\\
 Now in the following  we give an example nonreflexive Arnse regular Banach algebra which $\tilde{{Z}_1^\ell}(A^{**})= A^{**}.$\\\\
\noindent{\it{\bf Example 2-2.}} Consider the algebra $c_0=(c_0,.)$ is the collection of all sequences of scalars that convergence to $0$, with the some vector space operations and norm as $\ell_\infty$. We shows that  $\tilde{{Z}_1^\ell}(c_0^{**})=c_0^{**}$.\\

\begin{proof} We know that $c_0^*=\ell^1$, $c_0^{**}=\ell^\infty$ and $c_0^{***}=({\ell^\infty})^*={\ell^1}^{**}=\ell^1\oplus c_0^\perp$ as a Banach $c_0-bimodule$. Take $a^{\prime\prime\prime}\in c_0^{***}$. Then $a^{\prime\prime\prime}=(a^\prime,t)$ where $a^\prime\in \ell^1$ and $t\in c_0^\perp$. By [6, 2.6.22], $a^{\prime\prime}a^\prime\in c_0^\bot$ for all $a^{\prime\prime}\in c_0^{**}=\ell^\infty$. Then we have
$$<ta^{\prime\prime},a^\prime>=<t,a^{\prime\prime}a^\prime>=0.$$
It follows that $ta^{\prime\prime}=0$.
So for all $b^{\prime\prime}\in c_0^{**}=\ell^\infty$, we have the following equalities
$$<a^{\prime\prime\prime},a^{\prime\prime}b^{\prime\prime}>=<a^\prime,a^{\prime\prime}b^{\prime\prime}>\oplus <t,a^{\prime\prime}b^{\prime\prime}>=<a^{\prime\prime}b^{\prime\prime},a^\prime>\oplus <ta^{\prime\prime},b^{\prime\prime}>$$$$=<a^{\prime\prime}b^{\prime\prime},a^\prime>.$$
Consequently $\tilde{{Z}_1^\ell}(c_0^{**})={Z}_1(c_0^{**})$. Since $c_0^{**}$ is Arens regular, $\tilde{{Z}_1^\ell}(c_0^{**})=c_0^{**}$.\\
\end{proof}

\noindent{\it{\bf Theorem 2-3.}} Suppose that $A$ is a Banach algebra and $B\subseteq A^{**}$.
\begin{enumerate}
\item  If $A^{***}B\subseteq A^*$, then $B\subseteq \tilde{{Z}_1^\ell}(A^{**})$.
 \item  If $BA^{***}\subseteq A^*$, then $B\subseteq \tilde{{Z}_1^r}(A^{**})$.
\end{enumerate}
\begin{proof} 1) Let $b\in B$ and $(a_\alpha^{\prime\prime})_\alpha\subseteq A^{**}$ such that
$a_\alpha^{\prime\prime}\stackrel{w^*} {\rightarrow}a^{\prime\prime}$. We shows that
$ba_\alpha^{\prime\prime}\stackrel{w} {\rightarrow}ba^{\prime\prime}$. Let $a^{\prime\prime\prime}\in A^{***}$. Since $A^{***}B\subseteq A^*$, $a^{\prime\prime\prime}b\in A^*$. Then we have
$$<a^{\prime\prime\prime},ba_\alpha^{\prime\prime}>=<a^{\prime\prime\prime}b,a_\alpha^{\prime\prime}>=
<a_\alpha^{\prime\prime},a^{\prime\prime\prime}b>\rightarrow
<a^{\prime\prime},a^{\prime\prime\prime}b>=<a^{\prime\prime\prime},ba^{\prime\prime}>.$$
2) Proof is similar to (1).\\\\
\end{proof}

\noindent{\it{\bf Corollary 2-4.}} Suppose that $A$ is a Banach algebra. Then
\begin{enumerate}
\item If $A^{***}A^{**}\subseteq A^*$, then $\tilde{{Z}_1^\ell}(A^{**})={{Z}_1}(A^{**})= A^{**}$.
  \item If $A^{**}A^{***}\subseteq A^*$, then $\tilde{{Z}_1^r}(A^{**})= A^{**}$.\\
\end{enumerate}

\noindent{\it{\bf Corollary 2-5.}} Suppose that $A$ is a Banach algebra and $A^{***}A\subseteq A^*$. If $A$ is the left strong Arens irregular, then $\tilde{{Z}_1^\ell}(A^{**})=A$.\\\\
\noindent{\it{\bf Lemma 2-6.}} Suppose that $A$ is a Banach algebra. Let $(x^{\prime\prime}_\alpha)_\alpha\subseteq A^{**}$ and $(y^\prime_\beta)_\beta\subseteq A^*$ have taken $weak^*$ limit in $A^{**}$ and $A^{***}$, respectively. Let
$$\lim_\beta \lim_\alpha<x^{\prime\prime}_\alpha,y^\prime_\beta>=\lim_\alpha \lim_\beta <x^{\prime\prime}_\alpha,y^\prime_\beta>,$$
then we have the following assertions.\\
\begin{enumerate}
\item~~ $~\tilde{{Z}_1^r}(A^{**})=A^{**}$.
\item  $\tilde{{Z}_1^\ell}(A^{**})={{Z}_1}(A^{**})$.
\end{enumerate}
\begin{proof} 1) Let $b^{\prime\prime}\in A^{**}$. For all $a^{\prime\prime\prime}\in A^{***}$ there is $(a^\prime_\beta)_\beta\subseteq A^*$ such that $a^\prime_\beta\stackrel{w^*} {\rightarrow}a^{\prime\prime\prime}$. Assume that  $(a^{\prime\prime}_\alpha)_\alpha\subseteq A^{**}$ such that $a^{\prime\prime}_\alpha\stackrel{w^*} {\rightarrow}a^{\prime\prime}$. Then we have
$$<a^{\prime\prime\prime},a^{\prime\prime}b^{\prime\prime}>=\lim_\beta<a^\prime_\beta,a^{\prime\prime}b^{\prime\prime}>
=\lim_\beta \lim_\alpha<a^\prime_\beta,a_\alpha^{\prime\prime}b^{\prime\prime}>= \lim_\alpha \lim_\beta <a^\prime_\beta,a_\alpha^{\prime\prime}b^{\prime\prime}>$$
$$=
\lim_\alpha<a^{\prime\prime\prime},a_\alpha^{\prime\prime}b^{\prime\prime}>.$$
It follows that the mapping $a^{\prime\prime}\rightarrow a^{\prime\prime}b^{\prime\prime}~is~weak^*-to-weak~continuous~for~all~b^{\prime\prime}\in A^{**}$, and so $\tilde{{Z}_1^r}(A^{**})=A^{**}$.
2) Proof is similar to (1).

\end{proof}
\noindent By conditions  Lemma 2-6, it is clear that if $A$ is Arens regular, then  we have $\tilde{{Z}_1^\ell}(A^{**})=\tilde{{Z}_1^r}(A^{**})=A^{**}$.\\

\noindent{\it{\bf Theorem 2-7.}} Let $A$ be a Banach algebra and let  $B\subseteq A^{**}$. Then we have $B\tilde{{Z}_1^\ell}(A^{**})\subseteq \tilde{{Z}_1^\ell}(A^{**})$ and $\tilde{{Z}_1^r}(A^{**})B\subseteq \tilde{{Z}_1^r}(A^{**})$.

\begin{proof} We prove that $B\tilde{{Z}_1^\ell}(A^{**})\subseteq \tilde{{Z}_1^\ell}(A^{**})$ and proof of the other part is the same. Suppose that $a^{\prime\prime}\in  \tilde{{Z}_1^\ell}(A^{**})$ and $b^{\prime\prime}\in B$. Take
$a^{\prime\prime\prime}\in A^{***}$ and $(c_\alpha^{\prime\prime})_\alpha\subseteq A^{**}$ such that
$c^{\prime\prime}_\alpha\stackrel{w^*} {\rightarrow}c^{\prime\prime}$. Then
$$<a^{\prime\prime\prime},b^{\prime\prime}a^{\prime\prime}c^{\prime\prime}_\alpha>=
<a^{\prime\prime\prime}b^{\prime\prime},a^{\prime\prime}c^{\prime\prime}_\alpha>\rightarrow
<a^{\prime\prime\prime}b^{\prime\prime},a^{\prime\prime}c^{\prime\prime}>=
<a^{\prime\prime\prime},b^{\prime\prime}a^{\prime\prime}c^{\prime\prime}>.$$
It follows that the mapping $c^{\prime\prime}\rightarrow b^{\prime\prime}a^{\prime\prime}c^{\prime\prime}$ is $weak^*-to-weak$ continuous, and so $b^{\prime\prime}a^{\prime\prime}\in \tilde{{Z}_1^\ell}(A^{**})$.
\end{proof}

\noindent{\it{\bf Corollary 2-8.}} Let $A$ be a Banach algebra. Then we have the following assertions.
\begin{enumerate}
\item If $\tilde{{Z}_1^\ell}(A^{**})=A$, then $A$ is a right ideal in $A^{**}$.
\item If $\tilde{{Z}_1^r}(A^{**})=A$, then $A$ is a left ideal in $A^{**}$.
\item If $\tilde{{Z}_1^\ell}(A^{**})=\tilde{{Z}_1^r}(A^{**})=A$, then $A$ is an ideal in $A^{**}$.\\
\end{enumerate}

\noindent{\it{\bf Definition 2-9.}} Let $A$ be a Banach algebra. We define $\widetilde{wap}(A)$ as a subset of $A^{***}$ in the following
$$\widetilde{wap_\ell}(A)=\{a^{\prime\prime\prime}\in A^{***}:~a^{\prime\prime\prime}a^{\prime\prime}:A^{**}\rightarrow
\mathbb{C}~is ~weak^*~continuous~for~all~a^{\prime\prime}\in A^{**}\}.$$\\
It is clear that $\widetilde{wap_\ell}(A)=A^{***}$ if and only if $\tilde{{Z}_1^\ell}(A^{**})=A^{**}$. Thus, if  $\widetilde{wap_\ell}(A)=A^{***}$, then $wap(A)=A^*$, and so $A$ is Arens regular. It is easy to show that $wap(A)= \widetilde{wap_\ell}(A)$ if and only if $\tilde{{Z}_1^\ell}(A^{**})={{Z}_1}(A^{**})$.\\
The definition of $\widetilde{wap_r}(A)$ is similar.\\\\

\noindent{\it{\bf Theorem 2-10.}} Let $A$ be a Banach algebra. Then we have the following assertions.
\begin{enumerate}
\item If $A^{***}A^{**}\subseteq \widetilde{wap_\ell}(A)$, then $\tilde{{Z}_1^\ell}(A^{**})={{Z}_1}(A^{**})=A^{**}$.
\item If $A^{**}=A^{**}\tilde{{Z}_1^\ell}(A^{**})$, then $A^{***}A^{**}\subseteq \widetilde{wap_\ell}(A)$.
\end{enumerate}
\begin{proof} 1) Suppose that $a^{\prime\prime} \in A^{**}$ and $(b_\alpha^{\prime\prime})_\alpha\subseteq A^{**}$ such that
$b^{\prime\prime}_\alpha\stackrel{w^*} {\rightarrow}b^{\prime\prime}$. Take $a^{\prime\prime\prime}\in A^{***}$. Then, since $A^{***}A^{**}\subseteq \widetilde{wap_\ell}(A)$, we have
$$<a^{\prime\prime\prime},a^{\prime\prime}b^{\prime\prime}_\alpha>=
<a^{\prime\prime\prime}a^{\prime\prime},b^{\prime\prime}_\alpha>\rightarrow
<a^{\prime\prime\prime}a^{\prime\prime},b^{\prime\prime}>=
<a^{\prime\prime\prime},a^{\prime\prime}b^{\prime\prime}>.$$
It follows that $a^{\prime\prime}\in \tilde{{Z}_1^\ell}(A^{**})$, and so $A^{**}= \tilde{{Z}_1^\ell}(A^{**})$. Since $\tilde{{Z}_1^\ell}(A^{**})\subseteq {{Z}_1}(A^{**})$, the result is hold.\\
2) Let $a^{\prime\prime} \in A^{**}$ and $a^{\prime\prime\prime}\in A^{***}$. Since
$A^{**}=A^{**}\tilde{{Z}_1^\ell}(A^{**})$, there are $b^{\prime\prime} \in A^{**}$ and $c^{\prime\prime} \in\tilde{{Z}_1^\ell}(A^{**})$ such that $a^{\prime\prime}=b^{\prime\prime}c^{\prime\prime}$.
Consequently, we have
$$<a^{\prime\prime\prime}a^{\prime\prime},d^{\prime\prime}>=
<a^{\prime\prime\prime}b^{\prime\prime},c^{\prime\prime}d^{\prime\prime}>,$$
where $d^{\prime\prime} \in A^{**}$. Since $c^{\prime\prime} \in\tilde{{Z}_1^\ell}(A^{**})$,  $a^{\prime\prime\prime}a^{\prime\prime}\in \widetilde{wap_\ell}(A)$.

\end{proof}

\noindent{\it{\bf Corollary 2-11.}} Let $A$ be a Banach algebra. If $A^{**}=A^{**}\tilde{{Z}_1^\ell}(A^{**})$, then
$\tilde{{Z}^\ell_1}(A^{**})={{Z}_1}(A^{**})=A^{**}$.
\begin{proof} Since $wap(A)\subseteq \widetilde{wap_\ell}(A)$, by using Theorem 2.9, proof is hold.

\end{proof}

\noindent{\it{\bf Corollary 2-12.}} Let $A$ be a Banach algebra. Then,  if $A^{***}A^{**}\subseteq {wap}(A)$, then $\tilde{{Z}_1^\ell}(A^{**})={{Z}_1}(A^{**})=A^{**}$.\\

\noindent{\it{\bf Corollary 2-13.}} Let $A$ be a Banach algebra and $B$ be a subspace of $A^{**}$. If $B=B\tilde{{Z}_1^\ell}(B)$, then $\tilde{{Z}^\ell_1}(B)={{Z}_1}(B)=B$.\\\\

\begin{center}
\section{ \bf  Weak topological center of module actions  }
\end{center}
 \noindent In this section, we will study the definition of the weak topological centers on module actions and find some relations between it and reflexivity (spacial Arens regularity) of Banach algebras. For a Banach $A-bimodule$ $B$ we will study this definition on $n-th$ dual of $A$ and $B$. In some parts of this section, we have some conclusions in the dual groups, that is, for a locally compact   group  $G$ and a mixed unit $e^{\prime\prime}$ of  $L^1(G)^{**}$, we have
 ${Z}^\ell_{e^{\prime\prime}}(L^1(G)^{**})\neq L^1(G)^{**}$ and ${Z}^r_{e^{\prime\prime}}(L^1(G)^{**})\neq L^1(G)^{**}$ where ${Z}^\ell_{e^{\prime\prime}}(L^1(G)^{**})$ is a weak locally topological center of $L^1(G)^{**}$.\\
 We  have also conclusions as $ \tilde{{Z}^\ell}_{L^1(G)^{**}}(L^1(G)^{***})\subseteq L^1(G)^{*}$ and $ \tilde{{Z}^r}_{L^1(G)^{***}}(L^1(G)^{**})\subseteq L^1(G)$.\\
 On the other hand for Banach algebras $c_0$ and its second dual  $\ell^\infty$, we conclude that  $ \tilde{{Z}^\ell}_{\ell^\infty}((\ell^\infty)^{*})\subseteq \ell^\infty$ and
 $ \tilde{{Z}^r}_{(\ell^\infty)^{*}}(\ell^\infty)\subseteq c_0$.\\
 We also will extend some results from [8] into general situations, that is, for a locally compact   group  $G$ and  $n\geq 3$, we conclude the following equalities
\begin{enumerate}
\item  ${Z}^\ell_{{L^1(G)^{(n)}}}(L^1(G)^{**})= L^1(G)$, see [8].
\item  ${Z}^\ell_{{M(G)^{(n)}}}(L^1(G)^{**})= L^1(G)$.
\item  ${Z}^\ell_{{M(G)^{(n)}}}(M(G)^{**})= M(G)$.\\
\end{enumerate}

 \noindent Suppose that $A$ is a Banach algebra and $B$ is a Banach $A-bimodule$. According to [5, pp.27 and 28], $B^{**}$ is a Banach $A^{**}-bimodule$, where  $A^{**}$ is equipped with the first Arens product. So we can define the topological centers and weak topological centers  of module actions.\\
First for a Banach $A-bimodule$ $B$,  we  define the topological centers of the  left and right module actions of $A$ on $B$ as follows:

$${Z}^\ell_{A^{**}}(B^{**})={Z}(\pi_r)=\{b^{\prime\prime}\in B^{**}:~the~map~~a^{\prime\prime}\rightarrow \pi_r^{***}(b^{\prime\prime}, a^{\prime\prime})~:~A^{**}\rightarrow B^{**}$$$$~is~~~weak^*-weak^*~continuous\}$$
$${Z}^\ell_{B^{**}}(A^{**})={Z}(\pi_\ell)=\{a^{\prime\prime}\in A^{**}:~the~map~~b^{\prime\prime}\rightarrow \pi_\ell^{***}(a^{\prime\prime}, b^{\prime\prime})~:~B^{**}\rightarrow B^{**}$$$$~is~~~weak^*-weak^*~continuous\}$$
$${Z}^r_{A^{**}}(B^{**})={Z}(\pi_\ell^t)=\{b^{\prime\prime}\in B^{**}:~the~map~~a^{\prime\prime}\rightarrow \pi_\ell^{t***}(b^{\prime\prime}, a^{\prime\prime})~:~A^{**}\rightarrow B^{**}$$$$~is~~~weak^*-weak^*~continuous\}$$
$${Z}^r_{B^{**}}(A^{**})={Z}(\pi_r^t)=\{a^{\prime\prime}\in A^{**}:~the~map~~b^{\prime\prime}\rightarrow \pi_r^{t***}(a^{\prime\prime}, b^{\prime\prime})~:~B^{**}\rightarrow B^{**}$$$$~is~~~weak^*-weak^*~continuous\}.$$

 \noindent Now we define the weak topological centers of the right and left module actions of $A$ on $B$ as follows:\\

\noindent{\it{\bf Definition 3-1.}} Let $B$ a Banach $A-bimodule$. Then we define the weak topological centers of left and right module actions as in the following.\\
$$\tilde{{Z}}^\ell_{A^{**}}(B^{**})=\{ b^{\prime\prime}\in B:~a^{\prime\prime}\rightarrow b^{\prime\prime}a^{\prime\prime}~is~weak^*-weak~continuous~for~all~a^{\prime\prime}\in A^{**}\},$$
$$\tilde{ {Z}}^r_{A^{**}}(B^{**})=\{ b^{\prime\prime}\in B:~a^{\prime\prime}\rightarrow a^{\prime\prime}b^{\prime\prime}~is~weak^*-weak~continuous~for~all~a^{\prime\prime}\in A^{**}\},$$
$$\tilde{{Z}}^\ell_{B^{**}}(A^{**})=\{ a^{\prime\prime}\in B:~a^{\prime\prime}\rightarrow b^{\prime\prime}a^{\prime\prime}~is~weak^*-weak~~continuous~for~all~b^{\prime\prime}\in B^{**}\},$$
$$\tilde{ {Z}}^r_{B^{**}}(A^{**})=\{ a^{\prime\prime}\in A:~a^{\prime\prime}\rightarrow a^{\prime\prime}b^{\prime\prime}~is~weak^*-weak~~continuous~for~all~b^{\prime\prime}\in B^{**}\}.$$\\
 \noindent Let $A^{(n)}$ and  $B^{(n)}$  be $n-th~dual$ of $B$ and $A$, respectively. By [25], $B^{(n)}$ is a Banach $A^{(n)}-bimodule$ and we may therefore define the topological centers and weak topological centers of the right and left module action of
$A^{(n)}$ on  $B^{(n)}$ as similar above. when $n$ is odd, we define  $B^{(n)}B^{(n-1)}$ as a subspace of $A^{(n)}$. that is, for all $b^{(n)}\in B^{(n)}$ and $b^{(n-1)}\in B^{(n-1)}$, we define
$$<b^{(n)}b^{(n-1)},a^{(n-1)}>=<b^{(n)},b^{(n-1)}a^{(n-1)}>;$$
and if $n=0$, we take $A^{(0)}=A$ and $B^{(0)}=B$.\\

\noindent{\it{\bf Theorem 3-2.}} Let $B$ a Banach $A-bimodule$. Then for $n\geq 2$, we have the following assertions.
\begin{enumerate}

\item ~~~${Z}^\ell_{A^{(n)}}(B^{(n+1)})=B^{(n+1)}$ if and only if ${\tilde{{Z}}^r}_{A^{(n)}}(B^{(n)})=B^{(n)}$.
\item   ${Z}^\ell_{B^{(n)}}(A^{(n+1)})=A^{(n+1)}$ if and only if ${\tilde{{Z}}^r}_{B^{(n)}}(A^{(n)})=A^{(n)}$.
\item ~  ${Z}^r_{A^{(n)}}(B^{(n+1)})=B^{(n+1)}$ if and only if $\tilde{{Z}}^\ell_{A^{(n)}}(B^{(n)})=B^{(n)}$.
\item ~  ${Z}^r_{B^{(n)}}(A^{(n+1)})=A^{(n+1)}$ if and only if $\tilde{{Z}}^\ell_{A^{(n)}}(B^{(n)})=B^{(n)}$.
\end{enumerate}
\begin{proof} 1) Suppose that ${{{Z}^\ell}}_{A^{(n)}}(B^{(n+1)})=B^{(n+1)}$ and $b^{(n)}\in B^{(n)}$. We show that the mapping $a^{(n)}\rightarrow a^{(n)}b^{(n)}$ is $ weak^*-to-weak$ continuous. Assume that  $(a_\alpha^{(n)})_\alpha\subseteq A^{(n)}$ such that $a_\alpha^{(n)} \stackrel{w^*} {\rightarrow} a^{(n)}$. Then for all $b^{(n+1)}\in B^{(n+1)}$, we have $b^{(n+1)}a_\alpha^{(n)} \stackrel{w^*} {\rightarrow} b^{(n+1)}a^{(n)}$. It follows that
$$<b^{(n+1)},a_\alpha^{(n)}b^{(n)}>=<b^{(n+1)}a_\alpha^{(n)},b^{(n)}>\rightarrow <b^{(n+1)}a^{(n)},b^{(n)}>$$$$=
<b^{(n+1)},a^{(n)}b^{(n)}>.$$
Thus we conclude that $b^{(n)}\in {\tilde{{Z}}^r}_{A^{(n)}}(B^{(n)})$.\\
The converse is similarly.\\
Proof of (2), (3), (4) is similar to (1).
\end{proof}

 \noindent In the following example, for some non-reflexive Banach algebras such as $A$, we show that it is maybe that  ${\tilde{{Z}}^r}_{A^{{**}}}(A^{{**}})$ is equal to $A^{{**}}$ or no.\\

\noindent{\it{\bf Example 3-3.}}  Let $A$ be non-reflexive Banach space and let $<f,x>=1$ for some $f\in A^*$ and $x\in A$. We define the product on $A$ by $ab=<f,b>a$. It is clear that   $A$ is a Banach algebra with this product and it has  right identity $x$, see [ 5 ]. By easy calculation, for all $a^\prime \in A^*$,  $a^{\prime\prime}\in A^{**}$ and $a^{\prime\prime\prime}\in A^{***}$, we have
$$a^\prime a=<a^\prime , a>f,$$
$$a^{\prime\prime} a^\prime =<a^{\prime\prime}, f>a^\prime ,$$
$$a^{\prime\prime\prime} a^{\prime\prime}  =<a^{\prime\prime}, a^{\prime\prime}><.,f>.$$
Therefore we have ${{{Z}^\ell}}_{A^{{**} }}(A^{{***}})\neq A^{{***}}$. So by Theorem 3-2, we have  ${\tilde{{Z}}^r}_{A^{{**}}}(A^{{**}})\neq A^{{**}}$.\\
Similarly, if we define the product on $A$ as $ab=<f,a>b$ for all $a, ~b\in A$, then we have ${{{Z}^\ell}}_{A^{{**} }}(A^{{***}})=A^{{***}}$. By using Theorem 3-2, it follows that ${\tilde{{Z}}^r}_{A^{{**}}}(A^{{**}})= A^{{**}}$.\\

\noindent{\it{\bf Theorem 3-4.}} Let $B$ be  a Banach $A-bimodule$. Then for $n\geq 2$, we have the following assertions.
\begin{enumerate}
\item  ~If $B^{(n)}B^{(n-1)}\subseteq A^{(n-2)}$, then ${Z}^\ell_{A^{(n-1)}}(B^{(n-1)})={\tilde{{Z}^\ell}}_{A^{(n-1)}}(B^{(n-1)})$.

\item  If $B^{(n)}A^{(n-1)}\subseteq B^{(n-2)}$ and $B^{(n-1)}$ has a left unite  $A^{(n-1)}-module$, then $B$ is reflexive.
\end{enumerate}
\begin{proof} 1) It is clear that ${\tilde{{Z}^\ell}}_{A^{(n-1)}}(B^{(n-1)})\subseteq {Z}^\ell_{A^{(n-1)}}(B^{(n-1)})$. We show that for all $b^{(n-1)}\in {Z}^\ell_{A^{(n-1)}}(B^{(n-1)})$, the mapping $a^{(n-1)}\rightarrow b^{(n-1)}a^{(n-1)}$ is $weak^*-to-weak$ continuous. Suppose that $(a_\alpha^{(n-1)})_\alpha\subseteq A^{(n-1)}$ and  $a_\alpha^{(n-1)} \stackrel{w^*} {\rightarrow} a^{(n-1)}$. Since $b^{(n-1)}\in {Z}^\ell_{A^{(n-1)}}(B^{(n-1)})$,  $b^{(n-1)}a_\alpha^{(n-1)} \stackrel{w^*} {\rightarrow} b^{(n-1)}a^{(n-1)}$. Take $b^{(n)}\in B^{(n)}$. Then, since
$B^{(n)}B^{(n-1)}\subseteq A^{(n-2)}$, we have
$$<b^{(n)},b^{(n-1)}a_\alpha^{(n-1)}>=<b^{(n)}b^{(n-1)},a_\alpha^{(n-1)}>=<a_\alpha^{(n-1)},b^{(n)}b^{(n-1)}>$$$$\rightarrow
<a^{(n-1)},b^{(n)}b^{(n-1)}>=<b^{(n)}b^{(n-1)},a^{(n-1)}>=<b^{(n)},b^{(n-1)}a^{(n-1)}>.$$

2) Let $e^{(n-1)}$ be a left unit $A^{(n-1)}-module$ for $B^{(n-1)}$ and let $(b_\alpha^{(n-1)})_\alpha\subseteq B^{(n-1)}$ such that $b_\alpha^{(n-1)} \stackrel{w^*} {\rightarrow} b^{(n-1)}$. Suppose that $b^{(n)}\in B^{(n)}$. Since $b^{(n)}e^{(n-1)}\in  B^{(n-2)}$, we have
$$ <b^{(n)},b_\alpha^{(n-1)}>=<b^{(n)},e^{(n-1)}b_\alpha^{(n-1)}>=<b^{(n)}e^{(n-1)},b_\alpha^{(n-1)}>$$$$=
<b_\alpha^{(n-1)},b^{(n)}e^{(n-1)}>\rightarrow <b^{(n-1)},b^{(n)}e^{(n-1)}>=<b^{(n-1)},b^{(n)}>.$$
Consequently, the weak topology and $weak^*$ topology on $B^{(n-1)}$ is coincide, and so $B^{(n-1)}$ is reflexive which implies that $B$ is reflexive.

\end{proof}

\noindent{\it{\bf Corollary 3-5.}} Let $B$ be a Banach $A-bimodule$. Then, for every $n\geq 3$, we have the following assertions.
\begin{enumerate}
\item ~~If $A^{(n)}A^{(n-1)}\subseteq A^{(n-2)}$, then ${Z}_1^\ell{(A^{(n-1)})})={\tilde{{Z}_1^\ell}}{(A^{(n-1)})}$.
\item If $A^{(n)}A^{(n-1)}\subseteq A^{(n-2)}$ and $A^{(n-3)}$ has a bounded left approximate identity $(=BLAI)$, then $A$ is reflexive.\\
\end{enumerate}

\noindent{\it{\bf Lemma 3-6.}} Let $B$ be a  Banach $A-bimodule$. Then for $n\geq 2$, we have
\begin{enumerate}
\item If  $B^{(n-2)}$ has a bounded right approximate identity ($=BRAI$) $A^{(n-2)}-module$, then $B^{(n)}$ has  right unit $A^{(n)}-module$.
\item  If  $B^{(n-1)}$ has a bounded $BRAI$ $A^{(n-1)}-module$, then $B^{(n)}$ has  left unit $A^{(n+1)}-module$.\\
\end{enumerate}
\begin{proof} 1) Let $(e_\alpha^{(n-2)})_\alpha\subseteq A^{(n-2)}$ be a $BRAI$ for $B^{(n-2)}$. By Goldstines theorem there is a right unit $e^{(n)}\in A^{(n)}$ such that it is $weak^*$ closure of $(e_\alpha^{(n-2)})_\alpha$. Without lose generality, let $e_\alpha^{(n-2)} \stackrel{w^*} {\rightarrow}e^{(n)}$ in $A^{(n)}$. Assume that $b^{(n-1)}\in B^{(n-1)}$ and $b^{(n-2)}\in B^{(n-2)}$. Then
$$<b^{(n-1)},b^{(n-2)}>=\lim_\alpha<b^{(n-1)},b^{(n-2)}e_\alpha^{(n-2)}>=\lim_\alpha<e_\alpha^{(n-2)},b^{(n-1)}b^{(n-2)}>$$
$$=
<e^{(n)},b^{(n-1)}b^{(n-2)}>.$$
Then $e^{(n)}b^{(n-1)}=b^{(n-1)}$. Now let $b^{(n)}\in B^{(n)}$. Then we have
$$<b^{(n)}e^{(n)},b^{(n-1)}>=<b^{(n)},e^{(n)}b^{(n-1)}>=<b^{(n)},b^{(n-1)}>.$$
We conclude that $b^{(n)}e^{(n)}=b^{(n)}$.\\
2) Proof is similar to (1).

\end{proof}

\noindent{\it{\bf Corollary 3-7.}} Let $B$ be a   Banach $A-bimodule$. Then
$B^{(n-2)}$ has a BAI $A^{(n-2)}-module$ if and only if $B^{(n)}$ has  a unit $A^{(n)}-module$.\\

\noindent{\it{\bf Theorem 3-8.}} Let $B$ be a   Banach $A-bimodule$. Then we have the following assertions.
\begin{enumerate}
\item~~ If  $B^{(n+1)}B^{(n)}= A^{(n+1)}$ and ${\tilde{{Z}^\ell}}_{A^{(n)}}(B^{(n)})=B^{(n)}$, then $A$ is reflexive.
\item  Let $e^{(n)}\in A^{(n)}$ be a left unit  $A^{(n)}-module$ for $B^{(n)}$ and $e^{(n)}\in {\tilde{{Z}^\ell}}_{B^{(n)}}(A^{(n)})$. Then $B$ is reflexive.
\end{enumerate}
\begin{proof} 1) Let $(a_\alpha^{(n)})_\alpha\subseteq A^{(n)}$ and  $a_\alpha^{(n)} \stackrel{w^*} {\rightarrow} a^{(n)}$. Let $a^{(n+1)}\in A^{(n+1)}$. Since $B^{(n+1)}B^{(n)}\subseteq A^{(n+1)}$, there are
$b^{(n+1)}\in B^{(n+1)}$ and $b^{(n)}\in B^{(n)}$ such that $a^{(n+1)}=b^{(n+1)}b^{(n)}$. Thus, we have
$$<a^{(n+1)},a_\alpha^{(n)}>=<b^{(n+1)}b^{(n)},a_\alpha^{(n)}>=<b^{(n+1)},b^{(n)}a_\alpha^{(n)}>\rightarrow
<b^{(n+1)},b^{(n)}a^{(n)}>$$$$=<a^{(n+1)},a^{(n)}>.$$
We conclude that $A^{(n)}$ is reflexive, and so $A$ is reflexive.\\
2) Let $(b_\alpha^{(n)})_\alpha\subseteq B^{(n)}$ and $b_\alpha^{(n)} \stackrel{w^*} {\rightarrow} b^{(n)}$. Let
$b^{(n+1)}\in B^{(n+1)}$. Then
$$<b^{(n+1)},b_\alpha^{(n)}>=<b^{(n+1)},e^{(n)}b_\alpha^{(n)}>\rightarrow <b^{(n+1)},e^{(n)}b^{(n)}>=
<b^{(n+1)},b^{(n)}>.$$
Thus $B$ is reflexive.

\end{proof}
\noindent{\it{\bf Corollary 3-9.}} Assume that $B$ is a   Banach $A-bimodule$ and $B$ is reflexive. If $B^*B=A$, then $A$ is reflexive.\\

\noindent{\it{\bf Corollary 3-10.}} Assume that $A$ is a Banach algebra and $A^{**}$ has a left unit such that $\tilde{{Z}_1^\ell}(A^{**})$ consisting of it. Then $A$ is reflexive.\\

\noindent For Banach  $A-bimodule$ $B$ and every members of $A^{**}$ or $B^{**}$, in the following we introduce locally weak topological centers of $A^{**}$ and $B^{**}$ and we find some useful relations of them with reflexivity or some other useful conclusions in Banach  $A-bimodule$.\\
  Let $B$ be a Banach  $A-bimodule$ and suppose that  $a^{\prime\prime}\in A^{**}$. We say that $a^{\prime\prime}\rightarrow \pi_\ell^{***}(b^{\prime\prime}, a^{\prime\prime})$~is~$weak^*-to-weak^*$~point~~continuous, if for every net $(a_\alpha^{\prime\prime})_\alpha\subseteq A^{**}$ such that
$a^{\prime\prime}_\alpha\stackrel{w^*} {\rightarrow}a^{\prime\prime}$, it follows that $a^{\prime\prime}_\alpha b^{\prime\prime}   \stackrel{w^*} {\rightarrow}a^{\prime\prime} b^{\prime\prime}$.\\
 Suppose that  $B$ is a   Banach $A-bimodule$. Assume that $a^{\prime\prime}\in A^{**}$. Then we define  the  locally topological center of $a^{\prime\prime}$ on $B^{**}$ as follows
$${{Z}}^\ell_{a^{\prime\prime}}(B^{**})=\{ b^{\prime\prime}\in B^{**}:~a^{\prime\prime}\rightarrow \pi_\ell^{***}(b^{\prime\prime}, a^{\prime\prime})~is~weak^*-weak^*~point~~continuous\},$$
$${ {Z}}^r_{a^{\prime\prime}}(B^{**})=\{ b^{\prime\prime}\in B^{**}:~b^{\prime\prime}\rightarrow \pi_r^{t***}(a^{\prime\prime}, b^{\prime\prime})~is~weak^*-weak^*~point~~continuous\}.$$
The definition of ${{Z}}^\ell_{b^{\prime\prime}}(A^{**})$ and ${ {Z}}^r_{b^{\prime\prime}}(A^{**})$ where
$b^{\prime\prime}\in B^{**}$ are similar.\\
It is clear that $$\bigcap_{a^{\prime\prime}\in A^{**}}{{Z}}^{\ell,r}_{a^{\prime\prime}}(B^{**})={{Z}}^{\ell,r}_{A^{**}}(B^{**}),$$
$$\bigcap_{b^{\prime\prime}\in B^{**}}{{Z}}^{\ell,r}_{b^{\prime\prime}}(A^{**})={{Z}}^{\ell,r}_{B^{**}}(A^{**}).$$\\

\noindent{\it{\bf Definition 3-11.}} Let $B$ be a   Banach $A-bimodule$. Assume that $a^{\prime\prime}\in A^{**}$. Then we define  the locally weak   topological center of $a^{\prime\prime}$ on $B^{**}$  in the following
$$\tilde{{Z}}_{a^{\prime\prime}}^\ell(B^{**})=\{ b^{\prime\prime}\in B^{**}:~a^{\prime\prime}\rightarrow b^{\prime\prime}a^{\prime\prime}~~is~~weak^*-weak~~point~~continuous\},$$
$$\tilde{ {Z}}^r_{a^{\prime\prime}}(B^{**})=\{ b^{\prime\prime}\in B^{**}:~a^{\prime\prime}\rightarrow a^{\prime\prime}b^{\prime\prime}~~is~~weak^*-weak~~point~~continuous\}.$$
The definition of $\tilde{{Z}}^\ell_{b^{\prime\prime}}(A^{**})$ and $\tilde{ {Z}}^r_{b^{\prime\prime}}(A^{**})$ where
$b^{\prime\prime}\in B^{**}$ are similar.\\

It is clear that $$\bigcap_{a^{\prime\prime}\in A^{**}}\tilde{{Z}}^{\ell,r}_{a^{\prime\prime}}(B^{**})=\tilde{{Z}}^{\ell,r}_{A^{**}}(B^{**}),$$
$$\bigcap_{b^{\prime\prime}\in B^{**}}\tilde{{Z}}^{\ell,r}_{b^{\prime\prime}}(A^{**})=\tilde{{Z}}^{\ell,r}_{B^{**}}(A^{**}).$$\\

\noindent{\it{\bf Theorem 3-12.}} Let $B$ be a   Banach $A-bimodule$ and let $A^{**}$ has a right unit $e^{\prime\prime}$ such that $\tilde{{Z}}^\ell_{e^{\prime\prime}}(B^{**})=B^{**}$. If  $B^*$  factors on the left with respect to $A$, then $B$ is reflexive.

\begin{proof}  Let $(e_{\alpha})_{\alpha}\subseteq A$ be a BAI for $A$ such that  $e_{\alpha} \stackrel{w^*} {\rightarrow}e^{\prime\prime}$. Since $B^*$  factors on the left with respect to $A$, for all $b^{\prime}\in B^{*}$ there are $a\in A$ and $x^\prime\in B^*$ such that $x^\prime a=b^\prime$. Then for all $b^{\prime\prime}\in B^{**}$ we have
$$<\pi_\ell^{***}(e^{\prime\prime},b^{\prime\prime}),b^\prime>
=<e^{\prime\prime},\pi_\ell^{**}(b^{\prime\prime},b^\prime)>=
\lim_\alpha<\pi_\ell^{**}(b^{\prime\prime},b^\prime),e_{\alpha}>$$
$$=\lim_\alpha<b^{\prime\prime},\pi_\ell^{*}(b^\prime,e_{\alpha})>
=\lim_\alpha<b^{\prime\prime},\pi_\ell^{*}(x^\prime a,e_{\alpha})>$$
$$=\lim_\alpha<b^{\prime\prime},\pi_\ell^{*}(x^\prime ,ae_{\alpha})>
=\lim_\alpha<\pi_\ell^{**}(b^{\prime\prime},x^\prime) ,ae_{\alpha}>$$$$=<\pi_\ell^{**}(b^{\prime\prime},x^\prime) ,a>
=<b^{\prime\prime},b^{\prime}>.$$
Thus $\pi^{***}_\ell(e^{\prime\prime},b^{\prime\prime})=b^{\prime\prime}$ consequently $B^{**}$ has left unit $A^{**}-module$.
Since $\tilde{{Z}}^\ell_{e^{\prime\prime}}(B^{**})=B^{**}$, it is clear that $B$ is reflexive.\\
\end{proof}

 \noindent Assume that $B$ is  a   Banach $A-bimodule$ and   $e^{\prime\prime}\in A^{**}$ is a  mixed unit for $A^{**}$. We say that $e^{\prime\prime}$ is left mixed unit for $B^{**}$, if $\pi_\ell^{***}(e^{\prime\prime}, b^{\prime\prime})=\pi_\ell^{t***t}(b^{\prime\prime},e^{\prime\prime})=b^{\prime\prime}$ for all $b^{\prime\prime}\in B^{**}$.\\
The definition of the right mixed unit for $B$ is similar. So, it is clear that $e^{(n)}\in A^{(n)}$ is a left (resp. right) mixed unit for $B^{(n)}$ if and only if ${{Z}}^\ell_{e^{(n)}}(B^{(n)})=B^{(n)}$ (resp. ${ {Z}}^r_{e^{(n)}}(B^{(n)})=B^{(n)}$) whenever $n$ is even.\\\\

\noindent{\it{\bf Theorem 3-13.}} Let $B$ be a   Banach $A-bimodule$. Then for every $n\geq 2$, we have the following assertions.\\
\begin{enumerate}
\item ~~~ If $e^{(n)}\in A^{(n)}$ is a right unit for $B^{(n)}$ and  ${{Z}}^\ell_{e^{(n)}}(B^{(n)})=B^{(n)}$, then $$A^{(n-2)}B^{(n-1)}= B^{(n-1)}.$$
\item ~ If $e^{(n)}\in A^{(n)}$ is a left unit for $B^{(n)}$ and  ${ {Z}}^r_{e^{(n)}}(B^{(n)})=B^{(n)}$, then $$B^{(n-1)}A^{(n-2)}= B^{(n-1)}.$$
\item  Suppose that  $e^{(n)}\in A^{(n)}$ is a right unit for $B^{(n)}$ and  $\tilde{{Z}}^\ell_{e^{(n)}}(B^{(n)})=B^{(n)}$. Then $A^{(n-2)}B^{(n-1)}= B^{(n-1)}$ and $weak^*$ closure of  $A^{(n-2)}B^{(n+1)}$~ is~ $B^{(n+1)}.$
\item Suppose that $e^{(n)}\in A^{(n)}$ is a left unit for $B^{(n)}$ and  $\tilde{ {Z}}^r_{e^{(n)}}(B^{(n)})=B^{(n)}$. Then $B^{(n-1)}A^{(n-2)}= B^{(n-1)}$ and $weak^*$ closure of  $B^{(n+1)}A^{(n-2)}$ ~is ~$B^{(n+1)}.$
\end{enumerate}
\begin{proof}
1) Since $e^{(n)}$ is a right unit for $B^{(n)}$, by using Theorem 3-6, there is a $BRAI$ as  $(e_\alpha^{(n-2)})_\alpha\subseteq A^{(n-2)}$ for $A^{(n-2)}$ such that $e_\alpha^{(n-2)} \stackrel{w^*} {\rightarrow}e^{(n)}$ in $A^{(n)}$. Let  $b^{(n)}\in B^{(n)}$. Since ${{Z}}^\ell_{e^{(n)}}(B^{(n)})=B^{(n)}$, $b^{(n)}e_\alpha^{(n-2)} \stackrel{w^*} {\rightarrow}b^{(n)}e^{(n)}=b^{(n)}$. Then for all $b^{(n-1)}\in B^{(n-1)}$, we have
$$<b^{(n)},e_\alpha^{(n-2)}b^{(n-1)}>=<b^{(n)}e_\alpha^{(n-2)},b^{(n-1)}>\rightarrow <b^{(n)}e^{(n)},b^{(n-1)}>=$$$$
<b^{(n)},b^{(n-1)}>.$$
It follows that  $e_\alpha^{(n-2)}b^{(n-1)}\stackrel{w} {\rightarrow}b^{(n-1)}$.  Consequently, by using Cohen factori{Z}ation Theorem, we have  $A^{(n-2)}B^{(n-1)}= B^{(n-1)}$.\\
2) Proof is similar to (1).\\
3) Since $\tilde{{Z}}^\ell_{e^{(n)}}(B^{(n)})\subseteq{{Z}}^\ell_{e^{(n)}}(B^{(n)})$, it is clear that $A^{(n-2)}B^{(n-1)}= B^{(n-1)}$. We prove that the $weak^*$ closure of $A^{(n-2)}B^{(n+1)}$ is $B^{(n+1)}$.
Since $e^{(n)}$ is a right unit for $B^{(n)}$, by using Theorem 3-6, there is $(e_\alpha^{(n-2)})_\alpha\subseteq A^{(n-2)}$ such that $e_\alpha^{(n-2)} \stackrel{w^*} {\rightarrow}e^{(n)}$ in $A^{(n)}$. Since $\tilde{ {Z}}^r_{e^{(n)}}(B^{(n)})=B^{(n)}$, $b^{(n)}e_\alpha^{(n-2)}\stackrel{w} {\rightarrow}b^{(n)}e^{(n)}=b^{(n)}$ for all $b^{(n)}\in B^{(n)}$. Then for every $b^{(n+1)}\in B^{(n+1)}$, we have
$$<e_\alpha^{(n-2)}b^{(n+1)},b^{(n)}>=<b^{(n+1)},b^{(n)}e_\alpha^{(n-2)}>\rightarrow <b^{(n+1)},b^{(n)}e^{(n)}>=<b^{(n+1)},b^{(n)}>.$$
It follow that $e_\alpha^{(n-2)}b^{(n+1)}\stackrel{w^*} {\rightarrow}b^{(n+1)}$ and so proof is hold.\\
4) Proof is similar to (3).

\end{proof}

\noindent{\it{\bf Corollary 3-14.}} Assume that $A$ is a Banach algebra. Then we have the following assertions.
\begin{enumerate}
\item  If ${Z}^\ell_{e^{\prime\prime}}(A^{**})=A^{**}$, then $A^*$ factors on the right where $e^{\prime\prime}$ is a right unite for $A^{**}$.
\item  If ${Z}^r_{e^{\prime\prime}}(A^{**})=A^{**}$, then $A^*$ factors on the left where $e^{\prime\prime}$ is a left unite for $A^{**}$.\\
\end{enumerate}
\noindent{\it{\bf Corollary 3-15.}} Assume that $A$ is a Banach algebra and has a $BAI$. If $A^*$ not factors on the one side, then $A$ is not Arens regular.

\begin{proof} Since ${Z}^\ell_1(A^{**})\subseteq {Z}^\ell_{e^{\prime\prime}}(A^{**})$, by Theorem 3-13, proof is hold.

\end{proof}

\noindent{\it{\bf Example 3-16.}} Let $G$ be a locally compact   group and $e^{\prime\prime}$ be a mixed unit for $L^1(G)^{**}$. Then we have ${Z}^\ell_{e^{\prime\prime}}(L^1(G)^{**})\neq L^1(G)^{**}$ and ${Z}^r_{e^{\prime\prime}}(L^1(G)^{**})\neq L^1(G)^{**}$.\\

\noindent{\it{\bf Theorem 3-17.}} Assume that  $A$ is a Banach algebra. Then for $n\geq 0$, we have the following assertions.
\begin{enumerate}
\item  If $A^{(n)}$ has a $BRAI$, then $ \tilde{{Z}^\ell}_{A^{(n+2)}}(A^{(n+3)})\subseteq A^{(n+1)}$. Moreover, \\if $A^{(n+3)}A^{(n+1)}\subseteq A^{(n+1)}$, then $ \tilde{{Z}^\ell}_{A^{(n+2)}}(A^{(n+3)})= A^{(n+1)}$.
\item  If $A^{(n)}$ has a $BLAI$, then $ \tilde{{Z}^r}_{A^{(n+3)}}(A^{(n+2)})\subseteq A^{(n)}$. Moreover, if $A^{(n+4)}A^{(n)}\subseteq A^{(n+1)}$, then $ \tilde{{Z}^r}_{A^{(n+3)}}(A^{(n+2)})= A^{(n)}$.
\end{enumerate}
\begin{proof} 1) Assume that $A^{(n)}$ has a $BRAI$ as  $(e_\alpha^{(n)})_\alpha \subseteq A^{(n)}$ such that
$e_\alpha^{(n)} \stackrel{w^*} {\rightarrow}e^{(n+2)}$ in $A^{(n+2)}$ where $e^{(n+2)}$ is a right unite in $A^{(n+2)}$. Let $a^{(n+3)}\in \tilde{{Z}^\ell}_{A^{(n+2)}}(A^{(n+3)})$ and let $(a_\alpha^{(n+2)})_\alpha \subseteq A^{(n+2)}$ such that $a_\alpha^{(n+2)} \stackrel{w^*} {\rightarrow}a^{(n+2)}$ in $A^{(n+2)}$. Then we have
$$<a^{(n+3)},a_\alpha^{(n+2)}>=<a^{(n+3)},a_\alpha^{(n+2)}e^{(n+2)}>=<a^{(n+3)}a_\alpha^{(n+2)},e^{(n+2)}>$$
$$=<e^{(n+2)},a^{(n+3)}a_\alpha^{(n+2)}>\rightarrow <e^{(n+2)},a^{(n+3)}a^{(n+2)}>=<a^{(n+3)},a^{(n+2)}>.$$
It follows that $a^{(n+3)}:A^{(n+3)}\rightarrow \mathbb{C}$ is $weak^*$ continuous, and so $a^{(n+3)}\in A^{(n+1)}$
Consequently, we have $ \tilde{{Z}^\ell}_{A^{(n+2)}}(A^{(n+3)})\subseteq A^{(n+1)}$.\\
Now let $A^{(n+3)}A^{(n+1)}\subseteq A^{(n+1)}$. We show that the mapping $a^{(n+2)}\rightarrow a^{(n+1)}a^{(n+2)}$ is $weak^*-to-weak$ continuous for all $a^{(n+1)}\in A^{(n+1)}$. Assume that $(a_\alpha^{(n+2)})_\alpha \subseteq A^{(n+2)}$ such that  $a_\alpha^{(n+2)} \stackrel{w^*} {\rightarrow}a^{(n+2)}$ in $A^{(n+2)}$. Let $a^{(n+3)}\in A^{(n+3)}$. Then we have
$$<a^{(n+3)},a^{(n+1)}a_\alpha^{(n+2)}>=<a^{(n+3)}a^{(n+1)},a_\alpha^{(n+2)}>=<a_\alpha^{(n+2)},a^{(n+3)}a^{(n+1)}>$$
$$\rightarrow
<a^{(n+2)},a^{(n+3)}a^{(n+1)}>=<a^{(n+3)},a^{(n+1)}a^{(n+2)}>.$$
It follows that $a^{(n+1)}\in \tilde{{Z}^\ell}_{A^{(n+2)}}(A^{(n+3)})$, and so $ \tilde{{Z}^\ell}_{A^{(n+2)}}(A^{(n+3)})= A^{(n+1)}$.\\
2) It is similar to (1).

\end{proof}

\noindent{\it{\bf Corollary 3-18.}} Assume that  $A$ is a Banach algebra. Then
\begin{enumerate}
\item  If $A$ has a $BRAI$ and $ \tilde{{Z}^\ell}_{A^{**}}(A^{{***}})= A^{{***}}$. Then $A$ is reflexive.
\item  If $A$ has a $BLAI$ and $ \tilde{{Z}^r}_{A^{***}}(A^{{**}})= A^{{**}}$. Then $A$ is reflexive.\\
\end{enumerate}
\noindent{\it{\bf Example 3-19.}} Let $G$ be a locally compact infinite  group.  Then  we have
\begin{enumerate}
\item  $ \tilde{{Z}^\ell}_{L^1(G)^{**}}(L^1(G)^{***})\subseteq L^1(G)^{*}$ and $ \tilde{{Z}^r}_{L^1(G)^{***}}(L^1(G)^{**})\subseteq L^1(G)$.
\item  Since $c^*_0=\ell^1$, we have $ \tilde{{Z}^\ell}_{\ell^\infty}((\ell^\infty)^{*})\subseteq \ell^\infty$ and
 $ \tilde{{Z}^r}_{(\ell^\infty)^{*}}(\ell^\infty)\subseteq c_0$.\\
\end{enumerate}
\noindent{\it{\bf Theorem 3-20.}} Assume that $B$ is  a   Banach $A-bimodule$. Then for all $n\geq 1$, we have
\begin{enumerate}
\item  If $B^{(n)}B^{(n-1)}= A^{(n-1)}$, then $ {Z}^\ell_{B^{(n)}}(A^{(n)})\subseteq {Z}_1(A^{(n)})$.
\item  If $B^{(n+2)}B^{(n+1)}= A^{(n+1)}$ and $B^{(n)}B^{(n+1)}= A^{(n-1)}$, then $ \tilde{{Z}}^r_{B^{(n+1)}}(A^{(n)})\subseteq \tilde{{Z}}^r_1(A^{(n)})$.\\
\end{enumerate}
\begin{proof} 1) Let $a^{(n)}\in   {Z}_{B^{(n)}}(A^{(n)})$. We show that $x^{(n)}\rightarrow a^{(n)}x^{(n)}$ from $A^{(n)}$ into $A^{(n)}$ is $weak^*-to-weak^*$ continuous. Suppose that $(x_\alpha^{(n)})_\alpha \subseteq A^{(n)}$ such that  $x_\alpha^{(n)} \stackrel{w^*} {\rightarrow}x^{(n)}$ in $A^{(n)}$. First we show that
$x_\alpha^{(n)}b^{(n)} \stackrel{w^*} {\rightarrow}x^{(n)}b^{(n)}$ in $B^{(n)}$. Let $b^{(n-1)} \in B^{(n-1)}$. Then $$<x_\alpha^{(n)}b^{(n)},b^{(n-1)}>=<x_\alpha^{(n)},b^{(n)}b^{(n-1)}>\rightarrow <x^{(n)},b^{(n)}b^{(n-1)}>=<x^{(n)}b^{(n)},b^{(n-1)}>.$$
Now let $a^{(n-1)}\in A^{(n-1)}$. Then we have
$$<a^{(n)}x_\alpha^{(n)},a^{(n-1)}>=<a^{(n)}x_\alpha^{(n)},b^{(n)}b^{(n-1)}>=<a^{(n)}x_\alpha^{(n)}b^{(n)},b^{(n-1)}>$$$$
\rightarrow <a^{(n)}x^{(n)}b^{(n)},b^{(n-1)}>=<a^{(n)}x^{(n)},b^{(n)}b^{(n-1)}>=<a^{(n)}x^{(n)},a^{(n-1)}>.$$
It follows that $a^{(n)}\in {Z}_1(A^{(n)})$.\\
2) Proof is similar to (1).\\
\end{proof}

\noindent{\it{\bf Corollary 3-21.}} Assume that $B$ is  a   Banach $A-bimodule$ with left and right module actions $\pi_\ell$ and $\pi_r$, respectively. Then
\begin{enumerate}
\item  If $\pi^{**}_\ell$  is surjective, then  ${Z}^\ell_{B^{**}}(A^{**})\subseteq {Z}_1(A^{**})$
\item  If $\pi^{t**}_r$  is surjective, then  ${Z}^r_{B^{**}}(A^{**})\subseteq {Z}_2(A^{**})$\\

\end{enumerate}

\noindent{\it{\bf Example 3-22.}} Let $G$ be a locally compact   group and $n\geq 3$. Then
\begin{enumerate}
\item  ${Z}^\ell_{{L^1(G)^{(n)}}}(L^1(G)^{**})= L^1(G)$, see [8].
\item  ${Z}^\ell_{{M(G)^{(n)}}}(L^1(G)^{**})= L^1(G)$.
\item  ${Z}^\ell_{{M(G)^{(n)}}}(M(G)^{**})= M(G)$.\\

\end{enumerate}

\noindent{\it{\bf Theorem 3-23.}} Let $B$ be a   Banach $A-bimodule$. Then for $n\geq 2$, we have the following assertions.\\
\begin{enumerate}
\item  ~If  $B^{(n-1)}$ has a  $BRAI$ $A^{(n-2)}-module$ such that $weak^*$ convergence to $e^{(n)}\in A^{(n)}$ and $ \tilde{{Z}}^r_{e^{(n)}}(B^{(n-1)})=B^{(n-1)}$, then $B^{(n)}$ has  right unit $A^{(n)}-module$.

   \item  ~If  $B^{(n-1)}$ has a  $BLAI$ $A^{(n-1)}-module$ such that $weak^*$ convergence to $e^{(n+1)}\in A^{(n+1)}$ and $ \tilde{{Z}}^\ell_{e^{(n+1)}}(B^{(n-1)})=B^{(n-1)}$, then $B^{(n)}$ has  right unit $A^{(n+1)}-module$.

\item  ~If  $B^{(n-2)}$ has a  $BLAI$ $A^{(n-2)}-module$  such that $weak^*$ convergence to $e^{(n)}\in A^{(n)}$ and  $ \tilde{{Z}}^\ell_{e^{(n)}}(B^{(n-1)})=B^{(n-1)}$, then $B^{(n)}$ has left unit $A^{(n)}-module$.

\item  ~Suppose that  $B^{(n-1)}$ has a  $BRAI$ $A^{(n-1)}-module$ such that $weak^*$ convergence to $e^{(n+1)}\in A^{(n+1)}$. Let  $ \tilde{{Z}}^\ell_{e^{(n+1)}}(B^{(n-1)})=B^{(n-1)}$. Then $e^{(n+1)}$ is a left unit $A^{(n+1)}-module$ for $B^{(n)}$.\\
\end{enumerate}
\begin{proof} 1) Suppose that $(e_\alpha^{(n-2)})_\alpha\subseteq A^{(n-2)}$  is a $BLAI$  for $B^{(n-1)}$. Then there is a right unit element  $e^{(n)} \in A^{(n)}$ for $A^{(n)}$ such that $e_\alpha^{(n-2)} \stackrel{w^*} {\rightarrow}e^{(n)}$ in $A^{(n)}$. Let $b^{(n-1)}\in B^{(n-1)}$.
Since $ \tilde{{Z}}^r_{e^{(n)}}(B^{(n-1)})=B^{(n-1)}$, we have $e_\alpha^{(n-2)}b^{(n-1)} \stackrel{w} {\rightarrow}e^{(n)}b^{(n-1)}$. Then for every $b^{(n)}\in B^{(n)}$, we have
$$<b^{(n)}e^{(n)},b^{(n-1)}>=<b^{(n)},e^{(n)}b^{(n-1)}>=\lim_\alpha<b^{(n)},e_\alpha^{(n-2)}b^{(n-1)}>$$$$
=<b^{(n)},b^{(n-1)}>.$$
The proof of (2), (3) and (4) is similar to (1).\\
\end{proof}

\noindent{\it{\bf Corollary 3-24.}} Let $B$ be a   Banach $A-bimodule$. Then for $n\geq 2$, we have the following assertions.\\
\begin{enumerate}
\item  ~If  $B^{(n-1)}$ has a  $BAI$ $A^{(n-2)}-module$ such that $weak^*$ convergence to $e^{(n)}\in A^{(n)}$ and $ \tilde{{Z}}^r_{e^{(n)}}(B^{(n-1)})=B^{(n-1)}$, then $B^{(n)}$ is unital $A^{(n)}-module$.
 \item  ~If  $B^{(n-1)}$ has a  $BAI$ $A^{(n-1)}-module$ such that $weak^*$ convergence to $e^{(n+1)}\in A^{(n+1)}$ and $ \tilde{{Z}}^\ell_{e^{(n+1)}}(B^{(n-1)})=B^{(n-1)}$, then $B^{(n)}$ is unital $A^{(n+1)}-module$.

\item  ~If  $B^{(n-2)}$ has a  $BAI$ $A^{(n-2)}-module$  such that $weak^*$ convergence to $e^{(n)}\in A^{(n)}$ and  $ \tilde{{Z}}^\ell_{e^{(n)}}(B^{(n-1)})=B^{(n-1)}$, then $B^{(n)}$ is unital $A^{(n)}-module$.

\end{enumerate}

\begin{proof} By using Lemma 3-6 and Theorem 3-23, proof is hold.

\end{proof}

\begin{center}
\section{ \bf Weak amenability of Banach algebras  }
\end{center}
In this section, for $n\geq 2$, we show that if the weak topological center of $A^{(n)}$ with respect to the first Arens product be itself, then $A^{(n-2)}$ is weakly amenable whenever $A^{(n)}$ is weakly amenable. As corollary we show that if $\widetilde{wap_\ell}(A^{(n)})\subseteq A^{(n+1)}$, then the weakly amenability of $A^{(n)}$ implies  the weakly amenability of $A^{(n-2)}$ where $n\geq 2$. For  a   Banach $A-bimodule$ $B$, we deal with the question of when the second adjiont $D^{\prime\prime}:A^{**}\rightarrow B^{***}$ of $D:A\rightarrow B^*$ is again a derivation. This problem has been studied in [5, 9, 15] and we try to give answer to this question, that is, if the left topological center of a Banach algebra $A$ is itself, then $D:A\rightarrow B^*$ is a derivation whenever $D^{\prime\prime}:A^{**}\rightarrow B^{**}$ is a derivation.\\
At the beginning, we introduce some elementary notions and definitions as follows:\\
Let $B$ a   Banach $A-bimodule$. A derivation from $A$ into $B$is a bounded linear mapping $D:A\rightarrow B$ such that $$D(xy)=xD(y)+D(x)y~~for~~all~~x,~y\in A.$$
The space of continuous derivations from $A$ into $B$ is denoted by ${Z}^1(A,X)$.\\
Easy example of derivations are the inner derivations, which are given for each $b\in B$ by
$$\delta_b(a)=ab-ba~~for~~all~~a\in A.$$
The space of inner derivations from $A$ into $B$ is denoted by $N^1(A,X)$.
The Banach algebra $A$ is said to be a amenable, when for every Banach $A-bimodule$ $B$, the inner derivations are only derivations existing from $A$ into $B^*$. It is clear that $A$ is amenable if and only if $H^1(A,X^*)={Z}^1(A,X^*)/ N^1(A,X^*)=\{0\}$.\\
A Banach algebra $A$ is said to be a amenable, if every derivation from $A$ into $A^*$ is inner. Similarly, $A$ is weakly amenable if and only if $H^1(A,A^*)={Z}^1(A,A^*)/ N^1(A,A^*)=\{0\}$.\\

\noindent{\it{\bf Theorem 4-1.}} Assume that $A$ is a Banach algebra and $\tilde{{Z}_1^\ell}(A^{(n)})=A^{(n)}$ where $n\geq 2$. If $A^{(n)}$ is weakly amenable, then $A^{(n-2)}$ is weakly amenable.\\

\begin{proof} Suppose that $D\in {{Z}}^1(A^{(n-2)},A^{(n-1)})$. First we show that $$D^{\prime\prime}\in {{Z}}^1(A^{(n)},A^{(n+1)}).$$ Let $a^{(n)},~ b^{(n)}\in A^{(n)}$  and let $(a_\alpha^{(n-2)})_\alpha,~ (b_\beta^{(n-2)})_\beta\subseteq A^{(n)}$ such that $a_\alpha^{(n-2)}\stackrel{w^*} {\rightarrow}a^{(n)}$ and
$b_\beta^{(n-2)}\stackrel{w^*} {\rightarrow}b^{(n)}$, respectively. Since ${{Z}_1^\ell}(A^{(n)})=A^{(n)}$ ~we have ~ $$\lim_\alpha \lim_\beta a_\alpha^{(n-2)}D(b_\beta^{(n-2)})=a^{(n)}D^{\prime\prime}( b^{(n)}).$$
It is also clear that
  $$\lim_\alpha \lim_\beta D(a_\alpha^{(n-2)}b_\beta^{(n-2)})=D^{\prime\prime}( a^{(n)})b^{(n)}.$$ Since $D$ is continuous, we conclude that
 $$D^{\prime\prime}( a^{(n)}b^{(n)})=\lim_\alpha \lim_\beta D(a_\alpha^{(n-2)}b_\beta^{(n-2)})=
\lim_\alpha \lim_\beta a_\alpha^{(n-2)}D(b_\beta^{(n-2)})+$$
$$\lim_\alpha \lim_\beta D(a_\alpha^{(n-2)})b_\beta^{(n-2)}=
 a^{(n)}D^{\prime\prime}(b^{(n)})+D^{\prime\prime}( a^{(n)})b^{(n)}.$$
Since $A^{(n)}$ is weakly amenable, $D^{\prime\prime}$ is inner. It follows  that $D^{\prime\prime}( a^{(n)})=
a^{(n)}a^{(n+1)}-a^{(n+1)}a^{(n)}$ for some $a^{(n+1)}\in A^{(n+1)}$. Take  $a^{(n-1)}=a^{(n+1)}\mid_{A^{(n+1)}}$ and
$a^{(n-2)}\in A^{(n-2)}$. Then $$D(a^{(n-2)})=D^{\prime\prime}(a^{(n-2)})=a^{(n-2)}a^{(n-1)}-a^{(n-1)}a^{(n-2)}=\delta_{a^{(n-1)}}(a^{(n-2)}).$$
Consequently, we have  $H^1(A^{(n-2)},A^{(n-1)})=0$, and so $A^{(n-2)}$ is weakly amenable.

\end{proof}

\noindent{\it{\bf Corollary 4-2.}} Let $A$ be a Banach algebra and let $\widetilde{wap_\ell}(A^{(n-1)})\subseteq A^{(n)}$  whenever $n\geq1$. If  $A^{(n)}$ is weakly amenable, then  $A^{(n-2)}$ is weakly amenable.

\begin{proof} Since $\widetilde{wap_\ell}(A^{(n-1)})\subseteq A^{(n)}$, $\tilde{{Z}_1^\ell}(A^{(n)})=A^{(n)}$. Then,  by Theorem
4-1, proof is hold.

\end{proof}

\noindent{\it{\bf Corollary 4-3.}} Let $A$ be a Banach algebra  and ${{{Z}^\ell}}_{A^{(n)}}(A^{(n+1)})=A^{(n+1)}$ where $n\geq 2$. If  ${A^{(n)}}$ is weakly amenable, then ${A^{(n-2)}}$ is weakly amenable.\\

\noindent{\it{\bf Corollary 4-4.}} Let $A$ be a Banach algebra and let $D:A^{(n-2)}\rightarrow A^{(n-1)}$ be a derivation. Then $D^{\prime\prime}:A^{(n)}\rightarrow A^{(n+1)}$ is a derivation, if $\tilde{{Z}_1^\ell}(A^{(n)})=A^{(n)}$ for every $n\geq 2$.\\

\noindent{\it{\bf Theorem 4-5.}} Let $A$ be a Banach algebra and let $B$ be a closed subalgebra of $A^{(n)}$ that is consisting of $A^{(n-2)}$. If $\tilde{{Z}_1^\ell}(B)=B$ and   $B$ is weakly amenable, then $A^{(n-2)}$ is weakly amenable.

\begin{proof} Suppose that $D:A^{(n-2)}\rightarrow A^{(n-1)}$ is a derivation and $p:A^{(n+1)}\rightarrow B^\prime$ is the restriction map, defined by $P(a^{(n+1)})=a^{(n+1)}\mid_B$ for every $a^{(n+1)}\in A^{(n+1)}$. Since $\tilde{{Z}^\ell}(B)=B$, $\bar{D}=PoD^{\prime\prime}\mid_B:B\rightarrow B^\prime$ is a derivation. Since $B$ is weakly amenable, there is $b^\prime\in B^\prime$ such that $\bar{D}=\delta_{b^\prime}$. We take $a^{(n-1)}=b^\prime\mid_{A^{(n-1)}}$, then $D=\bar{D}$ on $A^{(n-1)}$. Consequently, we have $D=\delta_{a^{(n-1)}}$.

\end{proof}

\noindent{\it{\bf Corollary 4-6.}} Let $A$ be a Banach algebra. If $A^{***}A^{**}\subseteq A^{*}$ and $A$ is weakly amenable, then $A^{**}$ is weakly amenable.

\begin{proof} By using Corollary 2-4 and Theorem 3-2, proof is hold.

\end{proof}

\noindent{\it{\bf Corollary 4-7.}} Let $A$ be a Banach algebra and let $\tilde{Z}_1^\ell(A^{(n)})$ be weakly amenable
whenever $n\geq 2$. Then  $A^{(n-2)}$ is weakly amenable.\\\\

\noindent{\it{\bf Theorem 4-8.}} Let $B$ be a   Banach $A-bimodule$ and $D:A^{(n)}\rightarrow B^{(n+1)}$ be a derivation where $n\geq 0$. If $\tilde{{Z}^\ell}_{A^{(n)}}(B^{(n)})=B^{(n)}$, then $D^{\prime\prime}:A^{(n+2)}\rightarrow B^{(n+3)}$ is a derivation.

\begin{proof}
Let $x^{(n+2)},~ y^{(n+2)}\in A^{(n+2)}$  and let $(x_\alpha^{(n)})_\alpha,~ (y_\beta^{(n)})_\beta\subseteq A^{(n)}$ such that $x_\alpha^{(n)}\stackrel{w^*} {\rightarrow}x^{(n+2)}$ and
$y_\beta^{(n)}\stackrel{w^*} {\rightarrow}y^{(n+2)}$. Then for all $b^{(n)}\in B^{(n)}$, we have
 $b^{(n)}x_\alpha^{(n)}\stackrel{w} {\rightarrow}b^{(n)}x^{(n+2)}$. Consequently, we have
$$<x_\alpha^{(n)}D^{\prime\prime}(y^{(n+2)}),b^{(n)}>=<D^{\prime\prime}(y^{(n+2)}),b^{(n)}x_\alpha^{(n)}>\rightarrow <D^{\prime\prime}(y^{(n+2)}),b^{(n)}x^{(n+2)}>$$$$=<x^{(n+2)}D^{\prime\prime}(y^{(n+2)}),b^{(n)}>.$$
Also we have the following equality
$$<D^{\prime\prime}(x^{(n+2)})y_\beta^{(n)},b^{(n)}>=<D^{\prime\prime}(x^{(n+2)}),y_\beta^{(n)}b^{(n)}>\rightarrow
<D^{\prime\prime}(x^{(n+2)}),y^{(n+2)}b^{(n)}>$$$$=<D^{\prime\prime}(x^{(n+2)})y^{(n+2)},b^{(n)}>.$$
Since $D$ is continuous, it follows that
$$D^{\prime\prime}(x^{(n+2)}y^{(n+2)})=\lim_\alpha \lim_\beta D(x_\alpha^{(n)}y_\beta^{(n)})=
\lim_\alpha \lim_\beta x_\alpha^{(n)}D(y_\beta^{(n)})+$$$$\lim_\alpha \lim_\beta D(x_\alpha^{(n)})y_\beta^{(n)}=
x^{(n+2)}D^{\prime\prime}(y^{(n+2)})+D^{\prime\prime}(x^{(n+2)})y^{(n+2)}$$

\end{proof}

\noindent{\it{\bf Corollary 4-9.}} Let $B$ be  a   Banach $A-bimodule$ and $\tilde{{Z}}^\ell_{A^{**}}(B^{{**}})=B^{{**}}$. Then, $H^1(A,B^*)={0}$ if and only if
$H^1(A^{**},B^{***})={0}$.\\\\

\noindent{\it{\bf Corollary 4-10.}} Let $B$ be a   Banach $A-bimodule$ and $\tilde{{Z}^\ell}_{A^{{**}}}(B^{{**}})=B^{{**}}$. Let $D:A\rightarrow B^*$ be a derivation. Then $D^{\prime\prime}(A^{**})B^{**}\subseteq A^*$.

\begin{proof} By using Theorem 4-8 and [15, Corollary 4-3], proof is hold.

\end{proof}

\noindent{\it{\bf Corollary 4-10.}} Let $B$ be a   Banach $A-bimodule$ and let $D:A\rightarrow B^*$ be a derivation such that $D$ is surjective. Then  ${{Z}^\ell}_{A^{{**}}}(B^{{**}})=\tilde{{Z}}_{A^{{**}}}
^\ell(B^{{**}})$.

\begin{proof} By using Theorem 3-3 and  [15, Corollary 4-3], proof is hold.

\end{proof}

\noindent{\it{\bf Problems:}} \\
i) Let $G$ be a locally compact infinite  group. Find the sets $\tilde{Z_1}^\ell(L^1(G)^{**})=?$ and
$\tilde{Z_1}^\ell(M(G)^{**})=?$\\
ii) Assume that $A$ and $B$ are Banach algebra. If $\tilde{{Z}^\ell}(A^{(n)})=A^{(n)}$ and $\tilde{{Z}^\ell}(B^{(n)})=B^{(n)}$,  what  can  say about $\tilde{{Z}^\ell}(A^{(n)}\otimes B^{(n)})$?\\
iii)  Assume that $A$ and $B$ are Banach algebras. If $\tilde{{Z}^\ell}(A^{(n)}\otimes B^{(n)})=A^{(n)}\otimes B^{(n)}$, what can say about $\tilde{{Z}^\ell}(A^{(n)})$ and  $\tilde{{Z}^\ell}(B^{(n)})$?\\
iv) Let $G$ be a locally compact   group and $n\geq 3$. By notice to Theorem 3-20 and example 3-22, what can say about to the following sets.

 $\tilde{{Z}}^\ell_{{L^1(G)^{(n)}}}(L^1(G)^{**}),~$
 $\tilde{{Z}}^\ell_{{M(G)^{(n)}}}(L^1(G)^{**}),~$
  $\tilde{{Z}}^\ell_{{M(G)^{(n)}}}(M(G)^{**})$.\\

\bibliographystyle{amsplain}

\noindent Department of Mathematics, University of Mohghegh Ardabili, Ardabil, Iran\\
{\it Email address:} haghnejad@aut.ac.ir\\\\

\end{document}
